\documentclass[11pt]{article}
\usepackage{amsmath}
\usepackage{mathrsfs}
\usepackage{amsthm}
\usepackage{amsfonts}
\usepackage{amssymb}
\usepackage{bm}
\usepackage{mathtools}
\usepackage{color}
\usepackage{tikz}

\usepackage{enumerate}
\usepackage{enumitem}
\setlist[enumerate,1]{label=(\arabic*),font=\textup,
leftmargin=7mm,labelsep=1.5mm,topsep=0mm,itemsep=-0.8mm}
\setlist[enumerate,2]{label=(\alph*).,font=\textup,
leftmargin=7mm,labelsep=1.5mm,topsep=-0.8mm,itemsep=-0.8mm}

\usepackage[pdfstartview=FitH,bookmarksopen=false,
colorlinks=true,linkcolor=blue,citecolor=blue]{hyperref}

\newtheorem{theorem}{Theorem}[section]

\newtheorem{lemma}{Lemma}[section]

\newtheorem{conjecture}[theorem]{Conjecture}

\title{\bf Spectral extremal graphs for disjoint cliques \thanks {Research was partially supported by the National
Nature Science Foundation of China (grant numbers 11871329, 11971298)}}
\author { Zhenyu Ni$^{1}$, Jing Wang$^{2}$,  Liying Kang$^{2}$\thanks{\em Corresponding author. Email address: lykang@shu.edu.cn (L. Kang), 1051466287@qq.com(Z. Ni), wj517062214@163.com(J. Wang)}\\
{\small $^{1}$Department of Mathematics, Hainan University,
Haikou 570228, P.R. China}\\
{\small $^{2}$Department of Mathematics, Shanghai University,
Shanghai 200444, P.R. China}}

\date{}
\begin{document}

\maketitle

\begin{abstract}

  The $kK_{r+1}$ is the union of $k$ disjoint copies of $(r+1)$-clique. Moon [Canad. J. Math. 20 (1968) 95--102] and    Simonovits [Theory of Graphs (Proc. colloq., Tihany, 1996)]
  independently showed that if $n$ is sufficiently large, then $K_{k-1}\vee T_{n-k+1,r}$ is the unique extremal graph for $kK_{r+1}$.
  In this paper, we consider the graph which has the maximum spectral radius among all graphs without $k$ disjoint cliques. We prove that if $G$ attains the maximum spectral radius over all $n$-vertex $kK_{r+1}$-free graphs for sufficiently large $n$, then $G$ is isomorphic to $K_{k-1}\vee T_{n-k+1,r}$.


\medskip

\noindent{\bf AMS (2000) subject classification:}  05C50; 05C35
\end{abstract}

\section{Introduction}

In this paper, we consider only simple and undirected graphs. For two vertex disjoint graphs $G,H$, the {\sl union} of graph $G$ and $H$ is the graph
$G\cup H$ with vertex set $V(G)\cup V(H)$ and edge set $E(G)\cup E(H)$. In particular, we write $kG$ the vertex-disjoint union of $k$ copies of $G$.
The {\sl join} of $G$ and $H$, denoted by $G\vee H$,  is the graph obtained from $G\cup H$ by adding edges joining every vertex of $G$ to every vertex of $H$.
For two graphs $G$ and $F$, $G$ is called $F$-{\sl free} if it does not contain a copy of $F$ as a subgraph. For a fixed graph $F$, the Tur\'{a}n type extremal problem is to determine the maximum number of edges among all $n$-vertex $F$-free graphs, where the maximum number of edges is called the {\sl Tur\'{a}n number}, denoted by $\mathrm{ex}(n,F)$. An $F$-free graph on $n$ vertices is called an {\sl extremal graph} for $F$ if it has $\mathrm{ex}(n,F)$ edges, and the set of all extremal graphs is denoted by $\mathrm{Ex}(n,F)$.

Let $K_r(n_1,\ldots,n_r)$ be the complete $r$-partite graph with classes of sizes $n_1 ,\ldots , n_r$.  If $\sum_{i=1}^{r}n_i$ $=n$ and $|n_i-n_j|\leq 1$ for any $1\leq i<j\leq r$, then $K_r(n_1,\ldots,n_r)$ is called an $r$-partite {\sl Tur\'{a}n graph}, 
 denoted by $T_{n,r}$.  The well-known Tur\'{a}n Theorem states that the extremal graph corresponding to Tur\'{a}n number $\mathrm{ex}(n,K_{r+1})$ is $T_{n,r}$, i.e. $\mathrm{ex}(n,K_{r+1})=|E(T_{n,r})|$.
There are lots of researches on Tur\'{a}n type extremal problems (such as \cite{Bollobas1978,ChenGould2003,ErdosFuredi1995,Turan1941}).
Simonovits \cite{Simonovits1968} and Moon \cite{JWMoon1968} showed that if $n$ is sufficiently large, then $K_{k-1}\vee T_{n-k+1,r}$ is the unique extremal graph for $kK_{r+1}$.

The following spectral version of the Tur\'{a}n type  problem was proposed in Nikiforov \cite{Nikiforov2010extremalspectra}: What is the maximum spectral radius of a graph $G$ on $n$ vertices without a subgraph isomorphic to a given graph $F$?
Researches of the spectral Tur\'{a}n type extremal problem have drawn increasingly extensive interest (for example, see \cite{Nikiforov2008,BabaiGuiduli2009,Nikiforov2010Zarankiewicz,YuanWangZhai2012,ZhaiWang2012,ZhaiWangFang2020}). 
Nikiforov \cite{Nikiforov07} showed that if $G$ is a $K_{r+1}$-free graph on $n$ vertices,
then $\rho (G)\le \rho (T_{n,r})$, with equality if and only if
$G=T_{n,r}$. Cioab\u{a} et al. \cite{CioabaFengTaitZhang} proved that the spectral extremal graphs for $F_k$ belong to $\mathrm{Ex}(n,F_k)$, where $F_k$ is the graph consisting of $k$ triangles which intersect in exactly one common vertex.
Naturally, Cioab\u{a} et al. \cite{CDT21}
raised the following conjecture.

\begin{conjecture}(\cite{CDT21}) \label{conj}
Let $F$ be any graph such that the graphs in $\mathrm{Ex}(n,F)$ are Tur\'{a}n graphs plus $O(1)$ edges. Then for sufficiently large $n$, a graph attaining the maximum spectral radius among all $F$-free graphs on $n$ vertices is a member of $\mathrm{Ex}(n,F)$.
\end{conjecture}

The results of Nikiforov \cite{Nikiforov07}, Cioab\u{a} et al. \cite{CioabaFengTaitZhang}, Li et al. \cite{Yongtao21} and  Desai et al. \cite{DKLNTW} tell us that  Conjecture \ref{conj} holds for $K_{r+1}$, $F_k$, $H_{s,k}$ and $F_{k,r}$, where $H_{s,k}$ is the graph defined by intersecting $s$ triangles and $k$ odd cycles of length at least $5$ in exactly one common vertex, and $F_{k,r}$ is the graph consisting of $k$ copies of $K_r$ which intersect in a single vertex. Recently, Wang et al. \cite{WangKangXue}  proved Conjecture \ref{conj} completely.

In this paper, we shall prove the following theorem.

\begin{theorem}
\label{main}
For $k\geq 2$, $r\geq 2$, and sufficiently large $n$. Suppose that $G$ has the maximum spectral radius among all $kK_{r+1}$-free graphs on $n$ vertices, then $G$ is isomorphic to  $K_{k-1}\vee T_{n-k+1,r}$.
\end{theorem}


\section{Preliminaries}
\label{sec2}
Let $G=(V(G),E(G))$ be a simple connected graph with vertex set $V(G)$ and edge set $E(G)$.
For a vertex $v\in V(G)$, $N(v)$ is the set of neighbors of $v$ in $G$.  The {\sl degree} $d(v)$ of $v$ is $|N(v)|$, and the minimum and maximum degrees are denoted by $\delta(G)$ and $\Delta(G)$, respectively. We denote by $e(G)$ the number of edges  in $G$. For $V_1,V_2 \subseteq V(G)$, $E(V_1,V_2)$ denotes the set of edges of $G$ between  $V_1$ and  $V_2$, and $e(V_1,V_2)=|E(V_1,V_2)|$. For any $S\subseteq V(G)$, we write $N(S)=\cup_{u\in S}N(u)$, $d_{S}(v)=|N_{S}(v)|=|N(v)\cap S|$. Denote by $G\setminus S$ the graph obtained from $G$ by deleting all vertices in $S$ and their incident edges. $G[S]$ denotes the graph  induced by $S$ whose vertex set is $S$ and whose edge set consists of all edges of $G$ which have both ends in $S$. 
 A set $M$ of disjoint edges of $G$ is called a {\sl matching} in $G$.  The {\sl matching number}, denoted by $\nu(G)$, is the maximum cardinality of a matching in $G$. We call a matching with $k$ edges a {\sl $k$-matching}, denoted by $M_k$.
 For a matching $M$ of $G$, each vertex incident with an edge of $M$ is said to be {\sl covered} by $M$.

The {\sl adjacent matrix} of $G$ is $A(G)=(a_{ij})_{n\times n}$ with $a_{ij}=1$ if $ij\in E(G)$, and $a_{ij}=0$ otherwise.
The {\sl spectral radius} of $G$ is the largest eigenvalue of $A(G)$, denoted by $\rho(G)$.
For a connected graph $G$ on $n$ vertices, let $\mathbf{x}=(x_1,\ldots,x_n)^{\mathrm{T}}$ be an eigenvector of $A(G)$ corresponding to $\rho(G)$. Then $\mathbf{x}$ is a positive real vector, and
\begin{equation}\label{eigenequation}
\rho(G)x_i=\sum_{ij\in E(G)}x_j, \text{ for any } i\in [n].
\end{equation}
Another useful result concerns the Rayleigh quotient:
\begin{equation}\label{Rayleigh}
\rho(G)=\max_{\mathbf{x}\in \mathbb{R}^{n}_{+}}\frac{\mathbf{x}^{\mathrm{T}}A(G)\mathbf{x}}{\mathbf{x}^{\mathrm{T}}\mathbf{x}}=\max_{\mathbf{x}\in \mathbb{R}^{n}_{+}}\frac{2\sum_{ij\in E(G)}x_ix_j}{\mathbf{x}^{\mathrm{T}}\mathbf{x}}.
\end{equation}

The following spectral version of Stability Theorem was given  by Nikiforov \cite{Niki09JGT}.

\begin{theorem}[\cite{Niki09JGT}]  \label{lemniki}
Let $r\ge 2, 1/\ln n < c < r^{-8(r+21)(r+1)}, 0< \varepsilon < 2^{-36}r^{-24}$
and $G$ be a graph on $n$ vertices. If $\rho (G) > (1- \frac{1}{r} - \varepsilon )n$, then one of the following statements holds: \\
(a) $G$ contains a $K_{r+1}(\lfloor c\ln n\rfloor , \ldots ,\lfloor c\ln n\rfloor,
\lceil n^{1-\sqrt{c}}\rceil)$; \\
(b) $G$ differs from $T_{n,r}$ in fewer than $(\varepsilon^{1/4} +
c^{1/(8r+8)})n^2$ edges.
\end{theorem}

From the above theorem,  we can get the following result.
\begin{lemma}[\cite{DKLNTW}] \label{stability}
Let $F$ be a graph with chromatic number $\chi (F)=r+1$. For every $\varepsilon >0$, there exist $\delta >0$ and $n_0$ such that if  $G$ is an $F$-free graph on $n\ge n_0$ vertices  with $\rho (G) \ge (1- \frac{1}{r} -\delta )n$, then $G$ can be obtained from $T_{n,r}$ by adding and deleting at most $\varepsilon n^2$ edges.
\end{lemma}

Let $G$ be a simple graph with matching number $\nu(G)$ and maximum degree $\Delta(G)$. For two given integers $\nu$ and $\Delta$, define $f(\nu, \Delta)
=\max\{e(G): \nu(G)\leq \nu, \Delta(G)\leq \Delta \}$. In 1976, Chv\'atal and Hanson \cite{Chvatal76} obtained the following result.

\begin{lemma}[\cite{Chvatal76}]\label{ffunction}
For every two integers $\nu \geq 1$ and $\Delta \geq 1$, we have
$$f(\nu, \Delta)= \Delta \nu +\left\lfloor\frac{\Delta}{2}\right\rfloor
 \left \lfloor \frac{\nu}{\lceil{\Delta}/{2}\rceil }\right \rfloor
 \leq \Delta \nu+\nu.$$
\end{lemma}




The following lemma was given in \cite{CioabaFengTaitZhang}.
\begin{lemma}[\cite{CioabaFengTaitZhang}]
\label{intersect}
Let $V_1,\ldots,V_n$ be $n$ finite sets. Then
\[
|V_1 \cap \cdots \cap V_n| \geq \sum_{i=1}^{n}|V_i|-(n-1)|\bigcup_{i=1}^{n}V_i|
\]
\end{lemma}

\section{Proof of Theorem \ref{main}}
In this section we shall give a  proof of Theorem \ref{main}. Suppose that $G$ has the maximum spectral radius among all $kK_{r+1}$-free graphs on $n$ vertices, then we will  prove $G$ is isomorphic to $K_{k-1}\vee T_{n-k+1,r}$ for  sufficiently large $n$.
Clearly, $G$ is connected.
Let $\rho(G)$ be the spectral radius of $G$, $\mathbf{x}$ be a positive eigenvector of $\rho(G)$ with $\max\{x_i: i\in V(G)\}=1$. Without loss of generality, we assume $x_z=1$.

\begin{lemma}\label{rho}
Let $G$ be a $kK_{r+1}$-free graph on $n$ vertices with maximum spectral radius. Then
$$\rho(G)\geq \frac{r-1}{r}n+\frac{2(k-1)}{r}-\frac{1}{n}\left(\frac{(k-1)(r+k-1)}{r}+\frac{r}{2}\right).$$
\end{lemma}
\noindent{\bfseries Proof.} Let $H=K_{k-1}\vee T_{n-k+1,r}$.
Since $K_{k-1}\vee T_{n-k+1,r}$ is the unique extremal graph for $kK_{r+1}$, then
\begin{align}
\mathrm{ex}(n,kK_{r+1})&= e(T_{n-k+1,r})+(k-1)(n-k+1)+\binom{k-1}{2}\nonumber\\[2mm]
&\geq e(T_{n,r})+\frac{k-1}{r}n-\frac{(k-1)(r+k-1)}{2r}-\frac{r}{8}.\label{exnumber}
\end{align}

According to   (\ref{Rayleigh}) and (\ref{exnumber}),
 we have
\begin{align*}
\rho(G)&\geq \rho(H)
\geq \frac{\mathbf{1}^{\mathrm{T}}A(H)\mathbf{1}}{\mathbf{1}^{\mathrm{T}}\mathbf{1}}
= \frac{2\mathrm{ex}(n,kK_{r+1})}{n}\\[2mm]
&\geq \frac{2}{n}\left(e(T_{n,r})+\frac{k-1}{r}n-\frac{(k-1)(r+k-1)}{2r}-\frac{r}{8}\right)\\
&\geq\frac{r-1}{r}n+\frac{2(k-1)}{r}-\frac{1}{n}\left(\frac{(k-1)(r+k-1)}{r}+\frac{r}{2}\right).
\end{align*}
\qed

\begin{lemma}\label{partition}
Let $G$ be a $kK_{r+1}$-free graph on $n$ vertices with maximum spectral radius.
For every $\varepsilon >0$, there is an integer $n_0$ such that if $n\geq n_0$, then
\[
e(G)\geq e(T_{n,r})-\varepsilon n^2.
\]
Furthermore, $G$ has a partition $V(G)=V_1 \cup \cdots \cup V_{r}$ such that the number of crossing edges of $G$ $($i.e. $\sum_{1\leq i<j\leq r}e(V_i,V_j)$$)$ attains the maximum, and
\[
\sum_{i=1}^{r}e(V_i)\leq  \varepsilon n^2,
\]
and for any $i\in [r]$
$$\frac{n}{r}-3\sqrt{\varepsilon}n< |V_i|< \frac{n}{r}+3\sqrt{\varepsilon}n.$$
\end{lemma}
\noindent{\bfseries Proof.}
Since $G$ is $kK_{r+1}$-free,   
 by Lemmas \ref{stability} and  \ref{rho}, for sufficiently large $n$, there exists a partition of $V(G) = U_1 \cup \cdots \cup U_{r}$ such that $e(G) \geq e(T_{n,r}) - \varepsilon n^2$, $\sum_{i=1}^{r} e(U_i) \leq \varepsilon n^2$,
and $\lfloor\frac{n}{r}\rfloor \leq |U_i| \leq \lceil\frac{n}{r}\rceil$ for each $i \in [r]$.
Therefore, $G$ has a partition $V(G) = V_1 \cup \ldots \cup V_{r}$ such that the number of crossing edges of $G$ attains the maximum, and
$$\sum_{i=1}^{r} e(V_i) \leq \sum_{i=1}^{r} e(U_i) \leq \varepsilon n^2.$$
Let $a=\max\left\{\left||V_j|-\frac{n}{r}\right|, j\in [r]\right\}$.  Without loss of generality, we may assume that $\left||V_1|-\frac{n}{r}\right|=a$. Then
\begin{eqnarray*}
e(G)&\leq& \sum_{1\leq i<j\leq r}|V_i||V_j|+\sum_{i=1}^{r}e(V_i)\nonumber\\[2mm]
&\leq & |V_1|(n-|V_1|)+ \sum_{2\leq i<j\leq r}|V_i||V_j| +\varepsilon n^2\nonumber\\
&=& |V_1|(n-|V_1|)+ \frac{1}{2}\Big((\sum_{j=2}^{r}|V_j|)^2-\sum_{j=2}^{r}|V_j|^2\Big) +\varepsilon n^2\nonumber\\[2mm]
&\leq& |V_1|(n-|V_1|)+\frac{1}{2}(n-|V_1|)^2-\frac{1}{2(r-1)}(n-|V_1|)^2+\varepsilon n^2\\[2mm]
&<&-\frac{r}{2(r-1)}a^2+\frac{r-1}{2r}n^2+\varepsilon n^2,
\end{eqnarray*}
where  the last second inequality  holds by H\"{o}lder's inequality, and  the last inequality holds since $\left||V_1|-\frac{n}{r}\right|=a$.
On the other hand, since $e(G) \geq e(T_{n,r}) - \varepsilon n^2$, we have
$$e(G) \geq e(T_{n,r}) - \varepsilon n^2\geq \frac{r-1}{2r}n^2-\frac{r}{8}-\varepsilon n^2 > \frac{r-1}{2r}n^2-2\varepsilon n^2.$$ Therefore,
$\frac{r}{2(r-1)}a^2<3\varepsilon n^2$, which implies that $a<\sqrt{\frac{6(r-1)\varepsilon}{r}n^2}<3\sqrt{\varepsilon}n$. The proof is completed.
\qed

\begin{lemma}\label{W}
Suppose $\varepsilon$ and $\theta$ are two sufficiently small constants with $\theta< \frac{1}{20kr^4(r+1)}$ and $\varepsilon\leq \theta^2$. Let
\[
W:=\cup_{i=1}^{r}\{v\in V_i: d_{V_i}(v)\geq 2 \theta n\}.
\]
Then $|W|\leq \theta n$.
\end{lemma}
\noindent{\bfseries Proof.}
For all $i\in [r]$, let $W_i=W\cap V_i$. Then
\begin{equation*}
2e(V_i)=\sum_{u\in V_i}d_{V_i}(u)\geq \sum_{u\in W_i}d_{V_i}(u)\geq 2|W_i|\theta n.
\end{equation*}
Combining with Lemma \ref{partition}, we have
\begin{equation*}
\varepsilon n^2\geq \sum_{i=1}^{r}e(V_i)\geq  |W|\theta n,
\end{equation*}
which implies that $|W|\leq \frac{\varepsilon n}{\theta}\leq \theta n $.
\qed

\begin{lemma}\label{L}
Suppose $\varepsilon_1$ is a sufficiently small constant with $\sqrt{\varepsilon} < \varepsilon_1\ll \theta$. Let
\[
L:=\{v\in V(G): d(v)\leq (1-\frac{1}{r}-\varepsilon_1)n\}.
\]
Then $|L|\leq \varepsilon_2 n$, where $\varepsilon_2\ll \varepsilon_1$ is a sufficiently small constant satisfying $\varepsilon-\varepsilon_1 \varepsilon_2+\frac{r-1}{2r}\varepsilon_2^2<0$.
\end{lemma}
\noindent{\bfseries Proof.}
Suppose to the contrary that $|L|> \varepsilon_2 n$, then there exists $L'\subseteq L$ with $|L'|=\lfloor \varepsilon_2 n \rfloor$. Therefore,
\begin{eqnarray*}
e(G\setminus L')&\geq& e(G)-\sum_{v\in L'}d(v)\\[2mm]
&\geq & e(T_{n,r})-\varepsilon n^2- \varepsilon_2 n (1-\frac{1}{r}-\varepsilon_1)n\\[2mm]
&= & e(T_{n,r})-\varepsilon n^2-\frac{r-1}{r} \varepsilon_2 n^2+\varepsilon_1 \varepsilon_2 n^2\\[2mm]
&> & \frac{r-1}{2r}(n-\lfloor \varepsilon_2 n \rfloor)^2+\frac{k-1}{r}(n-\lfloor \varepsilon_2 n \rfloor)-\frac{(k-1)(k+r-1)}{2r}\\[2mm]
&\geq  & e(T_{n',r})+\frac{(k-1)n'}{r}-\frac{(k-1)(k+r-1)}{2r}\\[2mm]
&=  & \mathrm{ex}(n',kK_{r+1}),
\end{eqnarray*}
where $n'=n-\lfloor \varepsilon_2 n \rfloor$. Since $e(G\setminus L')>\mathrm{ex}(n-|L'|,kK_{r+1})$, $G\setminus L'$ contains a $kK_{r+1}$ as subgraph. This contradicts the fact that $G$ is $kK_{r+1}$-free.
\qed

\begin{lemma}\label{VminusWL}
For any $i\in [r]$, if $uv$ is an edge of $G[V_i\setminus (W\cup L)]$, then $G$ has $k(r+1)$ copies of $K_{r+1}$ which have only one common edge $uv$.
\end{lemma}
\noindent{\bfseries Proof.} For any $i\in [r]$, and any vertex $w\in V_i\setminus (W\cup L)$, we have $d(w)>(1-\frac{1}{r}-\varepsilon_1)n$, $d_{V_i}(w)<2\theta n$. Then for any $j\in [r]$ and $j\neq i$,
\begin{eqnarray*}
 d_{V_j}(w)&\geq& d(w)-d_{V_i}(w)-(r-2)(\frac{n}{r}+3\sqrt{\varepsilon}n)\\
&>& (1-\frac{1}{r}-\varepsilon_1)n- 2\theta n- (r-2)(\frac{n}{r}+3\sqrt{\varepsilon}n)\\
&> &\frac{n}{r}-3(r-1)\theta n.
\end{eqnarray*}
Without loss of generality, let $uv$ be an edge of $G[V_1\setminus (W\cup L)]$. We consider the common neighbors of $u,v$ in $V_2\setminus (W\cup L)$.
Combining with Lemma \ref{intersect}, we have
\begin{eqnarray*}
& &| N_{V_2}(u)\cap N_{V_2}(v)\setminus (W\cup L)|\\[2mm]
&\geq & d_{V_2}(u)+d_{V_2}(v)-|V_2|-|W|-|L|\\[2mm]
&> & 2(\frac{n}{r}-3(r-1)\theta n)-(\frac{n}{r}+3\sqrt{\varepsilon}n)-\theta n- \varepsilon_2 n\\[2mm]
&> & \frac{n}{r}-6r\theta n\\
&> & k(r+1).
\end{eqnarray*}
So there exist $k(r+1)$ vertices $u_{2,1},\ldots, u_{2,k(r+1)}$ in $V_2\setminus (W\cup L)$ such that the subgraph induced by two partitions $\{u,v\}$ and $\{u_{2,1},\ldots, u_{2,k(r+1)}\}$
is a complete bipartite graph. For an integer $s$ with $2\leq s\leq r-1$, suppose that there are vertices $u_{s,1},\ldots,u_{s,k(r+1)}\in V_{s}\setminus (W\cup L)$ such that $\{u,v\},$ $\{u_{2,1},\ldots,u_{2,k(r+1)}\},$ $\ldots,$ $\{u_{s,1},\ldots,u_{s,k(r+1)}\}$ induce a complete $s$-partite subgraph. We next consider the common neighbors of the above $(s-1)k(r+1)+2$ vertices in $V_{s+1}\setminus (W\cup L)$. 
By Lemma \ref{intersect}, we have
\begin{eqnarray*}
& &|N_{V_{s+1}}(u)\cap N_{V_{s+1}}(v)\cap(\cap_{i\in [s]\setminus\{1\}, j\in [k(r+1)]} N_{V_{s+1}}(u_{i,j}))\setminus (W\cup L)|\\[2mm]
&\geq & d_{V_{s+1}}(u)+d_{V_{s+1}}(v)+\sum_{i=2}^{s}\sum_{j=1}^{k(r+1)}d_{V_{s+1}}(u_{i,j})-((s-1)k(r+1)+1)|V_{s+1}|-|W|-|L|\\[2mm]
&> & ((s-1)k(r+1)+2)(\frac{n}{r}-3(r-1)\theta n)-((s-1)k(r+1)+1)(\frac{n}{r}+3\sqrt{\varepsilon}n)-\theta n- \varepsilon_2 n\\[2mm]
&> & \frac{n}{r}-12skr(r+1)\theta n\\
&> & k(r+1).
\end{eqnarray*}
 Then we can find  $k(r+1)$ vertices    $u_{s+1,1},\ldots, u_{s+1,k(r+1)}$ $\in V_{s+1}\setminus (W\cup L)$, which together with  $\{u,v\},$ $\{u_{2,1},\ldots,u_{2,k(r+1)}\},$ $\ldots,$ $\{u_{s,1},\ldots,u_{s,k(r+1)}\}$  forms a complete $(s+1)$-partite subgraph in $G$. Therefore, for every $2\leq i\leq r$, there exist $k(r+1)$ vertices in $V_{i}\setminus (W\cup L)$ such that $\{u_{2,1},\ldots,u_{2,k(r+1)}\},$ $\ldots,$ $\{u_{r,1},\ldots,u_{r,k(r+1)}\}$ induce a complete $(r-1)$-partite subgraph in $G$, and $u,v$ are adjacent to all the above $k(r-1)(r+1)$ vertices. Hence $G$ has $k(r+1)$ copies of $K_{r+1}$ which have only one common edge $uv$.
\qed

\begin{lemma}\label{independent}
For each $i\in [r]$, there exists an independent set $I_i\subseteq V_i\setminus (W\cup L)$ such that $|I_i|\geq |V_i\setminus (W\cup L)|-2(k-1)$.
\end{lemma}
\noindent{\bfseries Proof.} We first claim that $G[V_i\setminus (W\cup L)]$ is $M_k$-free for any $i\in [r]$. Suppose to the contrary that there exists $i_0\in [r]$ such that $G[V_{i_0}\setminus (W\cup L)]$ contains a copy of $M_k$. 
Then we can find a $kK_{r+1}$ by Lemma \ref{VminusWL}, and this contradicts the fact that $G$ is $kK_{r+1}$-free.
For every $i\in [r]$, let $M^i$ be a maximum matching of $G[V_i\setminus (W\cup L)]$, and $B^i$ be the set of vertices covered by $M^i$. Since  $G[V_i\setminus (W\cup L)]$ is $M_k$-free, $|B^i|\leq 2(k-1)$. Therefore, there exists an independent set $I_i\subseteq V_i\setminus (W\cup L)$ by deleting all vertices of $B^i$, and $|I_i|\geq |V_i\setminus (W\cup L)|-2(k-1)$.
\qed

\begin{lemma}\label{degreeViminusWL}
For any $i\in [r]$ and  any  $v\in V_i\setminus (W\cup L)$, $d_{V_i\setminus (W\cup L)}(v)<k(r+1)$.
\end{lemma}
\noindent{\bfseries Proof.} We will prove this lemma by contradiction. Without loss of generality, suppose that there exists a vertex  $u\in V_1\setminus (W\cup L)$ such that $d_{V_1\setminus (W\cup L)}(u)\geq k(r+1)$. Let $G'$ be the graph with $V(G')=V(G)$ and $E(G')=E(G)\cup \{uw: uw\notin E(G)\}$.
It follows from $u\in V_1\setminus (W\cup L)$ that $E(G)\subset E(G')$.
By the maximum of $\rho(G)$, $G'$ contains $kK_{r+1}$, say $F_1$, as a subgraph. From the construction of $G'$, we see that $u\in V(F_1)$, and there is a $(k-1)K_{r+1}$, say $F_2$, in $F_1\setminus \{u\}$. Obviously, $F_2\subseteq G$. Thus $F_2$ is a $(k-1)K_{r+1}$ copy of $G$, and $u\notin V(F_2)$. Since $d_{V_1\setminus (W\cup L)}(u)\geq k(r+1)$, there exists a vertex $v\in N_{V_1\setminus (W\cup L)}(u)$ such that $v\notin V(F_2)$. Then we can find $k(r+1)$ copies of $K_{r+1}$ which have only one common edge $uv$ by Lemma \ref{VminusWL}. Thus, we can find a $K_{r+1}$, say $F_3$, such that $V(F_3)\cap V(F_2)=\emptyset$. Thus $F_2\cup F_3$ is a $kK_{r+1}$ copy of $G$, which contradicts the fact that $G$ is $kK_{r+1}$-free.
\qed

\begin{lemma}\label{Fkrr+1}
For any $u\in W\setminus L$, $G$ contains $k(r+1)$ copies of $K_{r+1}$ which intersect only in $u$.
\end{lemma}
\noindent{\bfseries Proof.}
For any $u\in W\setminus L$, without loss of generality, we may assume that $u\in V_1$. Combining with Lemmas \ref{W} and \ref{L}, we have $d(u)>(1-\frac{1}{r}-\varepsilon_1)n$, and
\begin{eqnarray*}
 d_{V_1\setminus (W\cup L)}(u)&\geq& d_{V_1}(u)-|W\cup L|\\
&\geq & 2\theta n-\theta n-\varepsilon_2 n\\
&> & k(r+1).
\end{eqnarray*}
Let $u_{1,1},\ldots, u_{1,k(r+1)}$ be the neighbors of $u$ in $V_1\setminus (W\cup L)$. Then for every $i\in [k(r+1)]$, we have $d(u_{1,i})>(1-\frac{1}{r}-\varepsilon_1)n$, $d_{V_1}(u_{1,i})<2\theta n$, and
\begin{eqnarray}
 d_{V_2}(u_{1,i})&\geq& d(u_{1,i})-d_{V_1}(u_{1,i})-(r-2)(\frac{n}{r}+3\sqrt{\varepsilon}n)\nonumber\\
&> &\frac{n}{r}-\varepsilon_1 n-2\theta n-3(r-2)\sqrt{\varepsilon}n\nonumber\\
&> & \frac{n}{r}-3(r-1)\theta n.\label{eq4}
\end{eqnarray}
Since $V(G)=V_1 \cup \cdots \cup V_{r}$ is the vertex partition that  maximizes the number of crossing edges of $G$, we have $d_{V_1}(u)\leq \frac{1}{r}d(u)$. Therefore
\begin{eqnarray}
 d_{V_2}(u)&\geq& d(u)-d_{V_1}(u)-(r-2)(\frac{n}{r}+3\sqrt{\varepsilon}n)\nonumber\\
&> &\frac{r-1}{r}(1-\frac{1}{r}-\varepsilon_1)n-(r-2)(\frac{n}{r}+3\sqrt{\varepsilon}n)\nonumber\\
&>& \frac{n}{r^2}-\varepsilon_1 n-3(r-2)\sqrt{\varepsilon}n\nonumber\\
&>& \frac{n}{r^2}-(3r+5)\varepsilon_1n.\label{eq5}
\end{eqnarray}

We consider the common neighbors of $u,u_{1,1},\ldots,u_{1,k(r+1)}$ in $V_2\setminus (W\cup L)$.
Combining with Lemma \ref{intersect}, we have
\begin{eqnarray*}
& &|N_{V_2}(u)\cap (\cap_{i\in[k(r+1)]} N_{V_2}(u_{1,i}))\setminus (W\cup L)|\\[2mm]
&\geq & d_{V_2}(u)+\sum_{i=1}^{k(r+1)}d_{V_2}(u_{1,i})-k(r+1)|V_2|-|W|-|L|\\[2mm]
&> &\frac{n}{r^2}-(3r+5)\varepsilon_1 n + k(r+1)(\frac{n}{r}-3(r-1)\theta n)-k(r+1)(\frac{n}{r}+3\sqrt{\varepsilon}n)-\theta n- \varepsilon_2 n\\[2mm]
&> & \frac{n}{r^2}-16kr(r+1)\theta n\\
&> & k(r+1).
\end{eqnarray*}
Let $u_{2,1},\ldots, u_{2,k(r+1)}$ be the common neighbors of $u,u_{1,1},\ldots,u_{1,k(r+1)}$ in $V_2\setminus (W\cup L)$. For an integer $2\leq s\leq r-1$, suppose that $u_{s,1},\ldots,u_{s,k(r+1)}$ are the common neighbors of $\{u,u_{i,1},\ldots,u_{i,k(r+1)}: 1\leq i\leq s-1\}$ in $ V_{s}\setminus (W\cup L)$. 
We next consider the common neighbors of $\{u,u_{i,1},\ldots,u_{i,k(r+1)}: 1\leq i\leq s\}$ in $V_{s+1}\setminus (W\cup L)$. Using the similar method as in the proof of  (\ref{eq4}) and (\ref{eq5}), for every $i\in [s]$ and $j\in [k(r+1)]$, we have
\begin{equation*}
 d_{V_{s+1}}(u_{i,j})>\frac{n}{r}-3(r-1)\theta n,
\end{equation*}
and
\begin{equation*}
 d_{V_{s+1}}(u)>\frac{n}{r^2}-(3r+5)\varepsilon_1 n.
\end{equation*}
By Lemma \ref{intersect}, we have
\begin{eqnarray*}
& &|N_{V_{s+1}}(u)\cap (\cap_{i\in [s], j\in [k(r+1)]} N_{V_{s+1}}(u_{i,j}))\setminus (W\cup L)|\\[2mm]
&\geq & d_{V_{s+1}}(u)+\sum_{i=1}^{s}\sum_{j=1}^{k(r+1)}d_{V_{s+1}}(u_{i,j})-sk(r+1)|V_{s+1}|-|W|-|L|\\[2mm]
&> & \frac{n}{r^2}-(3r+5)\varepsilon_1n+ sk(r+1)(\frac{n}{r}-3(r-1)\theta n)-sk(r+1)(\frac{n}{r}+3\sqrt{\varepsilon} n)-\theta n- \varepsilon_2 n\\[2mm]
&> & \frac{n}{r^2}-16skr(r+1)\theta n\\
&> & k(r+1).
\end{eqnarray*}
Let $u_{s+1,1},\ldots, u_{s+1,k(r+1)}$ be the common neighbors of $\{u,u_{i,1},\ldots,u_{i,k(r+1)}: 1\leq i\leq s\}$ in $V_{s+1}\setminus (W\cup L)$. 
Therefore, for every $i\in [r]$, there exist $k(r+1)$ vertices, denoted by $\{u_{i,1},\ldots,u_{i,k(r+1)}\}$, in $V_{i}\setminus (W\cup L)$ such that $\{u_{1,1},\ldots,u_{1,k(r+1)}\}$, $\{u_{2,1},\ldots,u_{2,k(r+1)}\},$ $\ldots,$ $\{u_{r,1},\ldots,u_{r,k(r+1)}\}$ form a complete $r$-partite subgraph in $G$, and $u$ is adjacent to the above $kr(r+1)$ vertices. Hence we can find $k(r+1)$ copies of $K_{r+1}$ in $G$ which intersect only in $u$.
\qed

\begin{lemma}\label{WminusL}
$|W\setminus L|\leq k-1$.
\end{lemma}
\noindent{\bfseries Proof.}
Suppose to the contrary that $|W\setminus L|\geq k$. By Lemma \ref{Fkrr+1}, for any $u\in W\setminus L$, we can find $k(r+1)$ copies of $K_{r+1}$ in $G$ which intersect only in $u$.  Therefore, we can find at least $k$ disjoint $K_{r+1}$ in $G$. This is a contradiction  to the fact that $G$ is $kK_{r+1}$-free.
\qed

\begin{lemma}\label{Lemptyset}
$L=\emptyset.$
\end{lemma}
\noindent{\bfseries Proof.}
Let $x_{v_0}=\max\{x_v : v\in V(G)\setminus W\}$. Recall that $x_z=\max\{x_v : v\in V(G)\}=1$, then
\[
\rho(G)=\rho(G)x_z\leq |W|+(n-|W|)x_{v_0}.
\]
By Lemmas \ref{L} and \ref{WminusL}, we have
\begin{equation}\label{Wnumber}
|W|=|W\cap L|+|W\setminus L|\leq |L|+k-1\leq \varepsilon_2 n+k-1.
\end{equation}
Combining with Lemma \ref{rho}
, we have
\begin{eqnarray}\label{xv0}
x_{v_0}\geq \frac{\rho(G)-|W|}{n-|W|}
\geq \frac{\rho(G)-|W|}{n}
\geq 1-\frac{1}{r}- \varepsilon_2-\frac{O(1)}{n}
>1-\frac{2}{r}.
\end{eqnarray}
Therefore, we have
\begin{eqnarray*}
\rho (G) x_{v_0}&=&\sum_{vv_0\in E(G)} x_v
= \sum_{\substack{v\in W , vv_0\in E(G)}} x_v
+ \sum_{\substack{v\notin W  , vv_0\in E(G)} } x_v\\[2mm]
&\leq & |W|+(d(v_0)-|W|)x_{v_0},
\end{eqnarray*}
which implies that
\begin{eqnarray*}
d(v_0)&\geq& \rho(G)+|W|-\frac{|W|}{x_{v_0}}\\
&\geq& \rho(G)-\frac{2|W|}{r-2}\\
&\geq& \frac{r-1}{r}n+\frac{2(k-1)}{r}-\frac{1}{n}\left(\frac{(k-1)(r+k-1)}{r}+\frac{r}{2}\right)-\frac{2\varepsilon_2 n}{r-2}-\frac{2(k-1)}{r-2}\\
&>& (1-\frac{1}{r}-\varepsilon_1)n,
\end{eqnarray*}
where the last inequality  holds as $\varepsilon_2\ll \varepsilon_1$. Thus we have $v_0\notin L$, that is $v_0\in V(G)\setminus (W\cup L)$. Without loss of generality, we assume that $v_0\in V_1\setminus (W\cup L)$. Combining with Lemmas \ref{independent} and \ref{degreeViminusWL}, we have
\begin{eqnarray*}
 \rho (G) x_{v_0}  
 &=& \sum_{\substack{v\in W\cup L ,\\ vv_0\in E(G)}} x_v  + \sum_{\substack{v\in V_1\setminus (W\cup L) ,\\ vv_0\in E(G)}} x_v
 + \sum_{\substack{v\in (\cup_{i=2}^{r}V_i)\setminus (W\cup L),\\ vv_0\in E(G)} } x_v\\[2mm]
& < & |W|+|L|x_{v_0}+k(r+1)x_{v_0}+ \sum_{\substack{v\in \cup_{i=2}^{r}I_i , \\ vv_0\in E(G)} } x_v+ \sum_{\substack{v\in (\cup_{i=2}^{r}V_i\setminus I_i)\setminus (W\cup L)  ,\\ vv_0\in E(G)} } x_v\\[2mm]
& \leq & |W|+|L|x_{v_0}+k(r+1)x_{v_0}+ 2(k-1)(r-1)x_{v_0}+ \sum_{\substack{v\in \cup_{i=2}^{r}I_i } } x_v,
\end{eqnarray*}
which implies that
\begin{equation}\label{eigenvectorI}
\sum_{\substack{v\in \cup_{i=2}^{r}I_i } } x_v\geq (\rho (G)-|L|-k(3r-1)+2(r-1)) x_{v_0}-|W|.
\end{equation}

Next we will prove $L=\emptyset$. Suppose to the contrary that there is a vertex $u_0\in L$, then $d(u_0)\leq (1-\frac{1}{r}-\varepsilon_1)n$. 
Let $G'$ be the graph with $V(G')=V(G)$ and $E(G')=E(G\setminus \{u_0\}) \cup \{wu_0: w\in \cup_{i=2}^{r}I_i\}$. It is obvious that $G'$ is $kK_{r+1}$-free. Combining with Lemmas \ref{rho}, \ref{L}, (\ref{Wnumber}), (\ref{xv0}) and (\ref{eigenvectorI}), we have
\begin{align*}
\rho(G') - \rho(G)
&\geq \frac{\mathbf{x}^T\left(A(G')-A(G)\right)\mathbf{x}}{\mathbf{x}^T\mathbf{x}}
 = \frac{2x_{u_0}}{\mathbf{x}^T\mathbf{x}}\left(\sum_{\substack{w\in \cup_{i=2}^{r}I_i}} x_w - \sum_{uu_0\in E(G)} x_u\right) \\[2mm]
& \geq \frac{2x_{u_0}}{\mathbf{x}^T\mathbf{x}}\Bigl((\rho (G)-|L|-k(3r-1)+2(r-1)) x_{v_0}-2|W|-(d(u_0)-|W|)x_{v_0} \Bigr)\\[2mm]
& = \frac{2x_{u_0}}{\mathbf{x}^T\mathbf{x}}\Bigl((\rho (G)-|L|-k(3r-1)+2(r-1)-d(u_0)+|W|) x_{v_0}-2|W| \Bigr)\\[2mm]
&\geq \frac{2x_{u_0}}{\mathbf{x}^T\mathbf{x}}\Bigl( \frac{r-2}{r}(\varepsilon_1 n- \varepsilon_2 n- O(1))-2|W| \Bigr)\\[2mm]
&\geq \frac{2x_{u_0}}{\mathbf{x}^T\mathbf{x}}\Bigl( \frac{r-2}{r}(\varepsilon_1 n- \varepsilon_2 n- O(1))-2(\varepsilon_2 n+k-1) \Bigr)>0
\end{align*}
where the last inequality holds since $\varepsilon_2 \ll \varepsilon_1$. This contradicts the fact that $G$ has the largest spectral radius over all $kK_{r+1}$-free graphs, so $L$ must be empty.
\qed

\begin{lemma}\label{eigenvector}
For any $v\in V(G)$, $x_v\geq 1-\frac{1}{r-1}$.
\end{lemma}
\noindent{\bfseries Proof.}
Since $L=\emptyset$, then $|W|=|W\setminus L|\leq k-1$ by Lemma \ref{WminusL}. Let $x_{v_0}=\max\{x_v : v\in V(G)\setminus W\}$.
Recall that $x_z=\max\{x_v : v\in V(G)\}=1$, then
\[
\rho(G)=\rho(G)x_z\leq |W|+(n-|W|)x_{v_0}.
\]
Combining with Lemma \ref{rho}
, we have
\begin{eqnarray}\label{xv02}
x_{v_0}\geq \frac{\rho(G)-|W|}{n-|W|}
\geq \frac{\rho(G)-|W|}{n}
\geq 1-\frac{1}{r}-\frac{O(1)}{n}.
\end{eqnarray}
Using the similar method as in the proof of  (\ref{eigenvectorI}), we have
\begin{equation*}
\sum_{\substack{v\in \cup_{i=2}^{r}I_i } } x_v\geq (\rho (G)-k(r+3)+2) x_{v_0}-(k-1).
\end{equation*}
Suppose to the contrary that there exists $u\in V(G)$ such that $x_u< 1-\frac{1}{r-1}$. Let $G'$ be the graph with $V(G')=V(G)$ and $E(G')=E(G\setminus \{u\})\cup \{uw: w\in \cup_{i=2}^{r}I_i\}$. It is obvious that $G'$ is $kK_{r+1}$-free. Therefore, we have
\begin{align*}
\rho(G') - \rho(G)
&\geq \frac{\mathbf{x}^T\left(A(G')-A(G)\right)\mathbf{x}}{\mathbf{x}^T\mathbf{x}}
 = \frac{2x_{u}}{\mathbf{x}^T\mathbf{x}}\left(\sum_{\substack{w\in \cup_{i=2}^{r}I_i}} x_w - \sum_{uv\in E(G)} x_v\right) \\[2mm]
& \geq \frac{2x_{u}}{\mathbf{x}^T\mathbf{x}}\Bigl((\rho (G)-k(r+3)+2) x_{v_0}-(k-1)- \rho(G)x_u \Bigr)\\[2mm]
&> \frac{2x_{u}}{\mathbf{x}^T\mathbf{x}}\Bigl( (\rho (G)-k(r+3)+2)(1-\frac{1}{r}-\frac{O(1)}{n})-(k-1)- \rho(G)(1-\frac{1}{r-1}) \Bigr)\\[2mm]
&> \frac{2x_{u_0}}{\mathbf{x}^T\mathbf{x}}\Bigl( \frac{n}{r^2}-O(1) \Bigr)>0.
\end{align*}
This contradicts the fact that $G$ has the largest spectral radius over all $kK_{r+1}$-free graphs.
\qed

\begin{lemma}\label{Wk-1}
$|W|=k-1$, and $V_i\setminus W$ is an independent set for any $i\in [r]$.
\end{lemma}
\noindent{\bfseries Proof.}
Let $|W|=s$. Then $s\leq k-1$ by Lemmas \ref{WminusL} and \ref{Lemptyset}.

\noindent{\bfseries Claim } $\nu(\cup_{i=1}^{r}G[V_i\setminus W])\leq k-1-s$.

\noindent{\bfseries Proof of Claim.} Otherwise, $\nu(\cup_{i=1}^{r}G[V_i\setminus W])\geq k-s$. By Lemma \ref{VminusWL}, we can find a $(k-s)K_{r+1}$, denoted by $F_1$. Since $|W|=s$, by Lemma \ref{Fkrr+1}, we can find a $sK_{r+1}$, denoted by $F_2$, such that $V(F_1)\cap V(F_2)=\emptyset$. Therefore, $F_1\cup F_2$ is a  copy of $kK_{r+1}$ in  $G$,  a contradiction.

Suppose to the contrary that $s<k-1$.
By Lemmas \ref{degreeViminusWL} and \ref{Lemptyset}, we have
 $\Delta(\cup_{i=1}^{r}G[V_i\setminus W])<k(r+1)$. Combining with Lemma \ref{ffunction}, we have
\begin{align*}
e(\cup_{i=1}^{r}G[V_i\setminus W])&\leq f(\nu(\cup_{i=1}^{r}G[V_i\setminus W]),\Delta(\cup_{i=1}^{r}G[V_i\setminus W]))\\
&\leq f(k-s-1,k(r+1))\\
& \leq k(k-s)(r+1).
\end{align*}
Take $S\subseteq V_1\setminus W$ with $|S|=k-s-1$. Let $G'$ be the graph with $V(G')=V(G)$ and $E(G')=E(G)\setminus \{uv: uv\in \cup_{i=1}^{r}E(G[V_i\setminus W])\}\cup \{uv: u\in S,v\in (V_1\setminus W)\setminus S\}$. It is obvious that $G'$ is $kK_{r+1}$-free. 
Therefore,
\begin{align*}
&\rho(G') - \rho(G)\\[2mm]
&\geq \frac{\mathbf{x}^T\left(A(G')-A(G)\right)\mathbf{x}}{\mathbf{x}^T\mathbf{x}}\\[2mm]
&=\frac{2}{\mathbf{x}^T\mathbf{x}}\left(\sum_{ij\in E(G')}x_ix_j-\sum_{ij\in E(G)}x_ix_j\right) \\[2mm]
& \geq \frac{2}{\mathbf{x}^T\mathbf{x}}\left((k-s-1)(|V_1|-|W|-k+s+1)(1-\frac{1}{r-1})^2-k(k-s)(r+1)\right) \\[2mm]
& \geq \frac{2}{\mathbf{x}^T\mathbf{x}}\Bigl( (k-s-1)(\frac{n}{r}-3\sqrt{\varepsilon}n-k+1)(1-\frac{1}{r-1})^2-k(k-s)(r+1)\Bigr)\\[2mm]
&>0.
\end{align*}
This contradicts the fact that $G$ has the largest spectral radius over all $kK_{r+1}$-free graphs.
Therefore, $|W|=s=k-1$. Then it follows from the claim  that $\nu(\cup_{i=1}^{r}G[V_i\setminus W])\leq k-1-s=0$ for any $i\in [r]$.
So $V_i\setminus W$ is an independent set.
\qed

\begin{lemma}\label{WminusLdegree}
For any $u\in W$, $d(u)=n-1$.
\end{lemma}
\noindent{\bfseries Proof.}
Suppose to the contrary that there exists $u\in W$ such that $d(u)<n-1$. Let $v\in V(G)$ be a vertex such that $uv\notin E(G)$. Let $G'$ be the graph with $V(G')=V(G)$ and $E(G')=E(G)\cup \{uv\}$. We claim that $G'$ is $kK_{r+1}$-free. Otherwise, $G'$ contains a copy of $kK_{r+1}$, say $F_1$, as a subgraph, and $uv\in E(F_1)$. Let $F_2$ be the $K_{r+1}$ of $F_1$ which contains $uv$. 
Then $G$ contains a $(k-1)K_{r+1}$, denoted by $F_3$, as a subgraph, with $V(F_2)\cap V(F_3)=\emptyset$ and $u\notin V(F_3)$. Since $u\in  W$, by Lemmas \ref{Fkrr+1} and \ref{Lemptyset}, we can find a $K_{r+1}$, denoted by $F_4$, which contains $u$ and $V(F_4)\cap V(F_3)=\emptyset$. Thus $F_3\cup F_4$ is a copy of  $kK_{r+1}$ in $G$,  a contradiction. Therefore, $G'$ is $kK_{r+1}$-free. By the construction of $G'$, we have $\rho(G')>\rho(G)$, which contradicts the assumption that $G$ has the maximum spectral radius among all $kK_{r+1}$-free graphs on $n$ vertices.
\qed

\medskip
\noindent{\bfseries Proof of Theorem \ref{main}.}
Now we prove that $G$  is isomorphic to $ K_{k-1} \vee T_{n-k+1,r}$. For any $i\in[r]$, let $|V_i\setminus W|=n_i$.
By Lemmas \ref{Wk-1} and \ref{WminusLdegree}, there exists an $r$-partite graph $H$ with classes of size $n_1, n_2, \ldots, n_r$ such that
$G\cong K_{k-1} \vee H$.   By the maximum of $\rho(G)$, $H\cong K_{r}(n_1,n_2,\ldots,n_{r})$. It suffices to show that $|n_i-n_j|\leq 1$ for any $1\leq i<j\leq r$. Suppose  $n_1\geq n_2\geq \ldots \geq n_{r}$. We prove the assertion  by contradiction.  Assume that there exist $i_0, j_0$ with $1\leq i_0 < j_0\leq r$ such that $n_{i_0}-n_{j_0}\geq 2$. Let $H'=K_{r}(n_1,\ldots, n_{i_0}-1,\ldots,n_{j_0}+1,\ldots,n_{r})$, and $G'= K_{k-1} \vee H'$.

Recall that $\mathbf{x}$ is the eigenvector of $G$ corresponding to $\rho(G)$, by the symmetry we may assume $\mathbf{x}=(\underbrace{x_1,\ldots,x_1}_{n_1},\underbrace{x_2,\ldots,x_2}_{n_2},\ldots, \underbrace{x_r,\ldots,x_r}_{n_r},\underbrace{x_{r+1},\ldots,x_{r+1}}_{k-1})^{\mathrm{T}}$.
Thus by  (\ref{eigenequation}), we have
\begin{eqnarray}
\rho(G) x_i=\sum_{j=1}^{r}n_jx_j-n_ix_i+(k-1)x_{r+1}, \text{ for any } i\in [r], \label{14}
\end{eqnarray}
and \begin{eqnarray}
\rho(G) x_{r+1}=\sum_{j=1}^{r}n_jx_j+(k-2)x_{r+1}.\label{15}
\end{eqnarray}
Combining  (\ref{14}) and (\ref{15}), we have $x_i=\frac{\rho(G)+1}{\rho(G)+n_i}x_{r+1}$ for any $i\in [r]$, which implies that $x_{r+1}=\max\{x_v : v\in V(G)\}$. Recall that $\max\{x_v : v\in V(G)\}=1$,
 then $x_{r+1}=1$, and $x_i=\frac{\rho(G)+1}{\rho(G)+n_i}$ for any $i\in [r]$.
Let $u_{i_0}\in V_{i_0}\setminus W$ be a fixed vertex. Then $G'$ can be obtained from $G$ by deleting all edges between $u_{i_0}$ and $V_{j_0}\setminus W$, and adding all edges between $u_{i_0}$ and $V_{i_0}\setminus (W\cup\{u_{i_0}\})$.
According to  (\ref{Rayleigh}), we deduce that
\begin{align*}
\rho(G')-\rho(G)&\geq \frac{\mathbf{x}^{\mathrm{T}}(A(G')-A(G))\mathbf{x}}{\mathbf{x}^{\mathrm{T}}\mathbf{x}}\\[2mm]
&= \frac{2}{\mathbf{x}^{\mathrm{T}}\mathbf{x}}\left((n_{i_0}-1)x_{i_0}^2-n_{j_0}x_{i_0}x_{j_0}\right)\\[2mm]
&= \frac{2x_{i_0}}{\mathbf{x}^{\mathrm{T}}\mathbf{x}}\left((n_{i_0}-1)\frac{\rho(G)+1}{\rho(G)+n_{i_0}}-n_{j_0}\frac{\rho(G)+1}{\rho(G)+n_{j_0}}\right)\\[2mm]
&=\frac{2x_{i_0}}{\mathbf{x}^{\mathrm{T}}\mathbf{x}}\frac{(\rho(G)+1)(n_{i_0}\rho(G)-n_{j_0}\rho(G)-\rho(G)-n_{j_0})}{(\rho(G)+n_{i_0})(\rho(G)+n_{j_0})}\\[2mm]
&\geq\frac{2x_{i_0}}{\mathbf{x}^{\mathrm{T}}\mathbf{x}}\frac{(\rho(G)+1)(\rho(G)-n_{j_0})}{(\rho(G)+n_{i_0})(\rho(G)+n_{j_0})}
> 0,
\end{align*}
 where the last second inequality holds as $n_{i_0}-n_{j_0}\geq 2$, and the last inequality holds since $\rho(G)\geq \frac{r-1}{r}n+\frac{2(k-1)}{r}-\frac{1}{n}\left(\frac{(k-1)(r+k-1)}{r}+\frac{r}{2}\right)$, and $n_{j_0}=|V_{j_0}\setminus W|\leq \frac{n}{r}+3\sqrt{\varepsilon}n-(k-1)$. This contradicts the assumption that $G$ has the maximum spectral radius among all $n$-vertex $kK_{r+1}$-free graphs.
Therefore, $G$  is isomorphic to $ K_{k-1} \vee T_{n-k+1,r}$.
\qed

\end{document}